\documentclass[final]{siamltex}

\usepackage{epsfig}
\usepackage{tikz}
\usepackage{pgf}
\usepackage{url}
\usetikzlibrary{arrows,automata}
\usetikzlibrary{matrix,arrows,decorations.pathmorphing}
\usepackage[latin1]{inputenc}
\usepackage{tikz,fullpage}
\usepackage{tkz-berge}
\usepackage[position=top]{subfig}
\usepackage{verbatim}
\usetikzlibrary{arrows,shapes}
\usepackage{amssymb}
\usepackage{bm}

\def\e{{\epsilon}}

\def\0{{\bf 0}}
\def\1{{\bf 1}}
\def\e{{\bf e}}
\def\p{{\bf p}}

\def\+{{\bf +}}

\def\R{\mathbb{R}}

\def\Q{\mathbb{N}}

\providecommand{\com[2]}{\begin{tt}[#1: #2]\end{tt}}

\def\P{{\cal P}}

\newtheorem{remark}[theorem]{Remark}

\def\bang#1{\smallbreak\noindent$\triangleright$\ {\it #1}\ }

\def\ao#1{{#1}}
\def\aoc#1{{#1}}

\title{{\aoc How to decide consensus? a combinatorial necessary and sufficient condition and a proof that consensus is decidable but NP-hard}}

\author{Vincent D. Blondel\thanks{Division of Applied Mathematics, Universit\'e catholique de Louvain, 
4 avenue Georges Lemaitre, B-1348 Louvain-la-Neuve, Belgium. ({\tt vincent.blondel@uclouvain.be}).} \and Alex Olshevsky\thanks{Department of Industrial and Enterprise Systems Engineering, University of Illinois at Urbana Champaign, 104 S. Mathews Ave., Urbana, IL, 61801, USA ({\tt aolshev2@illinois.edu}). \newline A preliminary version of this paper was presented at the {\em 51st IEEE Conference on Decision and Control} in Dec. 2012 under the title {\em On the cost of deciding consensus}. Relative to the conference version, this version has complete proofs of all the theorems; additionally, the main theorem of this paper (Theorem \ref{necsufftheorem}) is new relative to the conference version.}}

\begin{document}

\maketitle

\begin{abstract} A set of stochastic matrices ${\cal P}$ is a consensus set if for every sequence of matrices $P(1), P(2), \ldots$ whose elements
belong to ${\cal P}$ and every initial state $x(0)$, the sequence of states defined by $x(t) = P(t) P(t-1) \cdots P(1) x(0)$ converges to a vector whose entries
are all identical. In this paper, we introduce an ``avoiding set condition" for compact sets of matrices and prove in our main theorem that this explicit combinatorial condition is both necessary and sufficient for consensus. We show that several of the conditions for consensus proposed in the literature can be directly derived from the avoiding set condition. The avoiding set condition is easy to check with an elementary algorithm, and so our result also establishes that  consensus  is algorithmically decidable. Direct verification of the avoiding set condition may require more than a polynomial time number of operations. This is however likely to be the case for any consensus checking algorithm since we also prove in this paper that unless $P=NP$,  consensus cannot be decided in polynomial time. 
\end{abstract}



\pagestyle{myheadings}
\thispagestyle{plain}

\section{Introduction} 

Recent years have seen considerable interest in the design of distributed control policies for networked systems of potentially  mobile agents. Motivated by
 the emerging applications of autonomous platoons of vehicles, formations of UAVs, coordinating mobile sensor platforms and other autonomous robotic systems,  the design of control strategies for swarms with unpredictable and time-varying connectivity has been the subject of much recent work. Controllers deployed in such systems ought to be completely decentralized, relying only on local information, and resilient to unexpected node and link failures. 

The design of such controllers usually relies on consensus algorithms, which are simple iterations allowing a distributed system of nodes to agree on a common
state with desirable properties. We mention recent work on coverage control \cite{gcb08}, formation control \cite{ofm07, othesis},
distributed estimation \cite{XBL05, XBL06}, distributed task assignment \cite{CBH09},  and distributed optimization \cite{TBA86} and \cite{NO09};
these papers and others design distributed control laws either by a direct reduction to an appropriately defined consensus problem or by using consensus 
protocols as a key subroutine.

Due to the wide applicability of consensus algorithms, the problem of finding convenient conditions that ensure their convergence has become an active area of current research. Much of the literature is focused on sufficient conditions for consensus algorithms to converge. Such conditions have become progressively more general with time; we refer the reader to \cite{degroot, cs77, TBA86, JLM03, M05, l05, ab06, lfm07, xw08, graphical, ll10, morse-ieee, charon, cao-sarym, cao-sarym2, julien-john} for a sampling.  As a result of these works, a number
of combinatorial properties of the underlying communication graphs are known to be sufficient for the convergence of classes of consensus algorithms. 

However, despite the significant attention attracted by this problem, no combinatorial necessary and sufficient conditions for consensus are known in general
(necessary and sufficient conditions are known in some special cases, for example when the matrices have positive diagonals \cite{morse-ieee}). Moreover,  much of the literature is formulated in terms of properties of matrix products which appear to be difficult to check,  so that it is unclear  whether it is possible to 
have efficiently verifiable conditions for consensus. Resolving these issues is the goal  of this paper. 

We next give an overview of our results, after setting down our notation and giving a quick motivating example. 

\subsection{Notation} Given a set of  stochastic matrices ${\cal P}$, an initial state $x(0)$, we consider the sequence of states $(x(0), x(1), \ldots)$ resulting from  
\begin{equation} \label{conscond} x(t)=P(t) x(t-1) \qquad t=1, 2, \ldots \end{equation}
where $P(1), P(2), \ldots$ are all matrices belonging to ${\cal P}$.  We will say that ${\cal P}$ is a \emph{consensus set} if such sequences $x(t)$ always converge to  a vector whose entries are all identical. In other words, being a consensus set requires that $\lim_{t \rightarrow \infty} x(t) = \alpha {\bf 1}$ for all matrix sequences $P(1), P(2), \ldots$ from ${\cal P}$ and initial states $x(0)$ (here ${\bf 1}$ denotes the vector of all ones {\ao and $\alpha$ is some scalar which may depend on $x(0)$ and on the matrix sequence}).  
As we will show in the next section, it is known that  ${\cal P}$ is a consensus set if and only if the limit 
\[ \lim_{t \rightarrow \infty} P(t) P(t-1) \cdots P(1)  \] 
always exists and has rank one.

\subsection{A motivating example} 

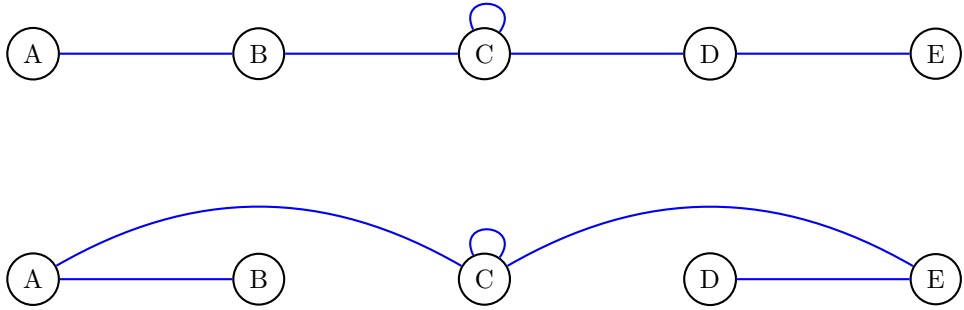
\begin{figure} 

\begin{center}
\begin{tikzpicture}[->, thick]
\SetVertexNormal[LineColor=black]
\SetVertexMath

\node (A1) at (0,3) [circle, draw] {A};
\node (B1) at (3,3) [circle, draw] {B};
\node (C1) at (6,3) [circle, draw] {C};
\node (D1) at (9,3) [circle, draw] {D};
\node (E1) at (12,3) [circle, draw] {E};

\path[-, color=blue]

(A1) edge (B1)
(B1) edge (C1)
(C1) edge (D1)
(D1) edge (E1)
(C1) edge [out=115, in=55, looseness=4] (C1);

\node (A2) at (0,0) [circle, draw] {A};
\node (B2) at (3,0) [circle, draw] {B};
\node (C2) at (6,0) [circle, draw] {C};
\node (D2) at (9,0) [circle, draw] {D};
\node (E2) at (12,0) [circle, draw] {E};

\path[-, color=blue]

(A2) edge (B2)
(A2) edge [bend left] (C2)
(C2) edge [bend left] (E2)
(D2) edge (E2)
(C2) edge [out=115, in=55, looseness=4] (C2);

\end{tikzpicture}
\caption{The graphs $G_1$ and $G_2$. The stochastic transition matrices corresponding to these two graphs do not form a consensus set.}
\label{gg}
\end{center}
 \end{figure}

Let us quickly illustrate with a small example that whether a set of matrices forms a consensus set may be subtle.  Consider the two graphs of Figure \ref{gg}. To a graph $G$ we associate a stochastic matrix $A$ by defining $a_{ij}=1/d(i)$ if $(i,j)$ is an edge in $G$, and $a_{ij}=0$ otherwise. Here $d(i)$ is the number of neighbors of node $i$, where we count node $i$ as its own neighbor if it has a self-loop.  To the  graphs of Figure \ref{gg} we associate to the two stochastic matrices $A_1, A_2$ and we claim that these do not form a consensus set. Indeed, it is easy to verify that the product $A_1 A_2$ has the property that the rows corresponding to nodes $A$ and $E$ have only a single nonzero entry, namely a $1$ on the diagonal. It therefore follows that the limit $\lim_k (A_1 A_2)^k$ does not approach a rank-$1$ matrix and the matrices $A_1$ and $A_2$ do not form a consensus set.

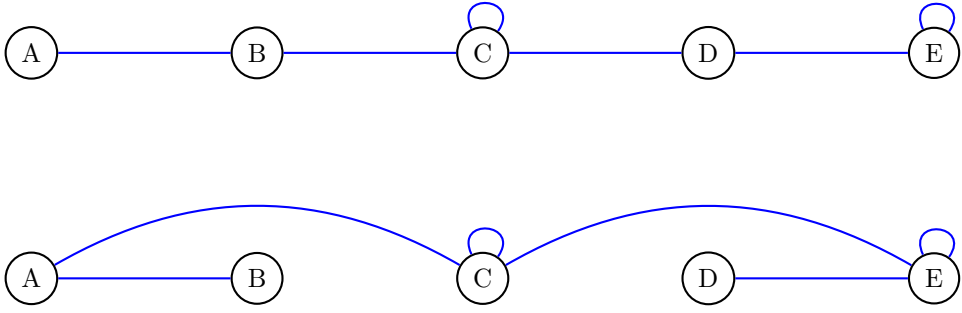
\begin{figure} \begin{center}

\begin{tikzpicture}[->, thick]
\SetVertexNormal[LineColor=black]
\SetVertexMath

\node (A3) at (0,-4) [circle, draw] {A};
\node (B3) at (3,-4) [circle, draw] {B};
\node (C3) at (6,-4) [circle, draw] {C};
\node (D3) at (9,-4) [circle, draw] {D};
\node (E3) at (12,-4) [circle, draw] {E};

\path[-, color=blue]

(A3) edge (B3)
(B3) edge (C3)
(C3) edge (D3)
(D3) edge (E3)
(C3) edge [out=115, in=55, looseness=4] (C3)
(E3) edge [out=115, in=55, looseness=4] (E3);

\node (A4) at (0,-7) [circle, draw] {A};
\node (B4) at (3,-7) [circle, draw] {B};
\node (C4) at (6,-7) [circle, draw] {C};
\node (D4) at (9,-7) [circle, draw] {D};
\node (E4) at (12,-7) [circle, draw] {E};

\path[-, color=blue]

(A4) edge (B4)
(A4) edge [bend left] (C4)
(C4) edge [bend left] (E4)
(D4) edge (E4)
(C4) edge [out=115, in=55, looseness=4] (C4)
(E4) edge [out=115, in=55, looseness=4] (E4);

\end{tikzpicture} 

\caption{The graphs $G'_1$ and $G'_2$. The only difference with the graphs on the previous figure is the additional self-loop at node $E$. The stochastic transition matrices corresponding to these two graphs form a consensus set.}
\label{ggg}
\end{center} \end{figure}

Consider now instead the  stochastic transition matrices associated to the two graphs $G'_1$ and $G'_2$ of Figure  \ref{ggg}. Note that the only difference with the graphs of Figure   \ref{gg} is the additional self-loop at node $E$. This turns out to make a difference for consensus: {\ao a simple computation reveals that every product of length $4$ of these two matrices has a column whose entries are strictly positive. It is well-known (and, in any case,  follows out of the more general results we will show - see Example 2 in Section \ref{examples})  that this implies these two matrices form a consensus set.}

\subsection{Our results}

In Section \ref{necsuff}, we present the ``avoiding sets condition'' which is a necessary and sufficient combinatorial condition for consensus. Roughly speaking, it says ${\cal P}$ is a consensus set if and only if one cannot find two cyclic set sequences which
avoid  each other in a finite graph sequence built from the matrices in ${\cal P}$. We refer the reader to Section \ref{necsuff} for a precise 
statement.  Moreover, we show that a number of the known results in the literature follow
immediately as corollaries of this fact. A new implication of this theorem is that checking whether a finite set of stochastic matrices is a consensus set is decidable.  

Unfortunately, checking whether a finite set of stochastic matrices is a consensus set by checking the avoiding sets condition is not  polynomial time.  In Section  \ref{hard}, we show that this is unavoidable. We first prove in Theorem  \ref{twotheorem} that, unless $P=NP$, the problem of determining if a given set of two stochastic matrices is a consensus set cannot be solved in polynomial-time. We then prove a similar result for undirected matrices, which
are nonnegative matrices $P$ with the property that if $p_{ij} > 0$ then $p_{ji} > 0$. Undirected matrices are important for consensus problems because they often appear in applications where information exchange goes in both directions: whenever $i$ communicates with $j$, $j$ communicates with $i$ (with possibly different strengths).
In Theorem \ref{threetheorem} we prove that the problem of determining if a finite set of stochastic undirected matrices is a consensus 
set cannot be decided in polynomial time (unless $P=NP$).
Our proof is for the case of three (or more) matrices. We believe that the same result holds true if there are only two undirected matrices but we have
not been able to derive the corresponding proof. We remark that this section is not independent but depends on the results of the previous sections.

Finally, we conclude our contribution with some remarks and open problems in Section \ref{concl}.
 
\section{The avoiding sets condition: a necessary and sufficient condition for consensus\label{necsuff}} The goal of this section is to prove a combinatorial necessary and sufficient condition for consensus which we call ``the avoiding sets condition,'' and discuss its implications. We begin with a quick definition before proceeding to state the theorem.

\smallskip

\begin{definition} Given a stochastic matrix $P \in \mathbb{R}^{n \times n}$ we define the directed graph $G(P)$ to have the vertex set $\{1, \ldots, n\}$ and the edge 
set $\{ (i,j) ~|~ p_{ij} > 0 \}$. If $S$ is a subset of $\{1, \ldots, n\}$, we will use $N_P(S)$ to denote the set of out-neighbors of $S$ in the graph $G(P)$. We will often be implicitly given a sequence of stochastic matrices $P(1), P(2), \ldots$, and in this case we will use $N_i(S)$ as shorthand for $N_{P(i)}(S)$. Moreover,  in this case $N^1(S)$ will be shorthand for $N_1(S)$, $N^2(S)$ will be shorthand for $N_2(N_1(S))$, $N^3(S)$ will be shorthand for $N_3(N_2(N_1(S)))$ and so forth.  
\end{definition}

\smallskip

With this notation in place, we can state the avoiding sets condition.

%

\bigskip

\begin{theorem}[The avoiding sets condition] \label{necsufftheorem} The compact set ${\cal P}$ of stochastic $n \times n$ matrices is not a consensus set if and only if there 
exist two sequences of nonempty subsets of $\{1, \ldots, n\}$,\[ S_1, S_2, \ldots, S_l, \] and \[ S_1', S_2', \ldots, S_l', \]  of length $l \leq 3^{n}-2^{n+1}+1$  and a sequence of 
matrices $P(1), P(2), \ldots, P(l)$ from ${\cal P}$ such that 
$$ S_i \cap S_i' = \emptyset ~~ \mbox{ for all } i = 1, \ldots, l,$$ and for all $k = 1, \ldots, l-1$,
  $$ N_{k}(S_k) \subset S_{k+1}, ~~~ N_{l}(S_l) \subset S_1$$
$$ N_{k}(S_k') \subset S_{k+1}', ~~~ N_{l}(S_l') \subset S_1'$$  
\end{theorem} 

\smallskip

\begin{remark} We summarize the essential meaning of the theorem in words. We will consider a finite list of $l$ sets to be a  ``set cycle'' if there is a corresponding sequence $P(1), \ldots, P(l)$ from ${\cal P}$ such that  the operation of replacing the $k$'th set with its out-neighbors in $G(P(k))$ maps it within the $k+1$'st set (and the last to the first). The theorem states that a set of matrices is a consensus set if and only if there do not exist two set-cycles which are disjoint at every step. 
\end{remark}

\smallskip

\begin{remark} Since finite sets of matrices are compact, we note that this theorem applies to all finite sets of stochastic matrices. 
\end{remark}

\smallskip

\begin{remark} Note that the two set-cycles have to be disjoint in the sense of $S_i \cap S_{i}' = \emptyset$. However, nothing prevents
$S_i$ and $S_j'$ from having nonzero intersection when $i \neq j$. Similarly, the two sets $S_i$ and $S_{j}$ are not required to have
empty intersection.
\end{remark}

\subsection{Examples and algorithmic implications\label{examples}} While Theorem \ref{necsufftheorem} may appear somewhat unwieldy at first glance, it can be used to establish
a variety of concrete results. We illustrate its use with several examples. 

\bigskip

\noindent {\bf Example 1:} We begin with an  elementary example. Suppose $P \in \R^{n \times n}$ is a ``disconnected'' stochastic matrix, meaning that we can partition  
$\{1, \ldots, n\}$ into two disjoint nonempty sets $U$ and $V$ such that \begin{equation} \label{noconn} p_{uv}=p_{vu}=0 ~~~~~~~~~~~\mbox{ for all } u \in U, v \in V 
\end{equation} In that case, 
$\{P\}$ is clearly not a consensus set. This obvious fact may be deduced from the avoiding sets condition by picking 
$S_1 = U, S_1' = V$; Eq. (\ref{noconn}) then insures that $N_P(S_1) \subset S_1, N_P(S_1') \subset S_1'$ as required.

%

\smallskip

Similarly, let us suppose $P \in \R^{n \times n}$ has the following bipartiteness-type property: there exist two disjoint nonempty sets $U_1, U_2$ of $\{1, \ldots, n\}$ such that \[ N_P(U_1) \subset U_2, N_P (U_2) \subset U_1. \]   It is easy to see that $\{P\}$ is not a consensus set,
and this may be deduced from the avoiding sets condition by choosing 
\[ S_1 = U_1, S_2 = U_2,\]
\[ S_1' = U_2, S_2'=U_1. \]

 
\bigskip

\noindent {\bf Example 2:} We continue with another elementary example. In the introduction, we mentioned the well-known fact that a finite set of matrices forms a consensus 
set if each has a strictly positive column. We justify that now. Indeed, if $P$ is a matrix with a positive column, then $N_P(U) \cap N_P(V) \neq 0$ for any two nonempty sets $U,V$, so  two avoiding set sequences cannot exist in this case. 

\smallskip

\noindent {\bf Example 3:} We outline how the classical consensus results \cite{TBA86, JLM03} may be obtained from the avoiding sets condition. First, consider the (compact) set  ${\cal P}_{\gamma}$ of stochastic matrices $P$ whose graphs $G(P)$ are strongly connected, which have positive diagonals, and whose every positive entry is furthermore 
at least $\gamma > 0$.  This is a consensus set, since regardless of the matrix sequence $P(1), P(2), \ldots$, we have that for any set $S$, $N^{n-1}(S) = \{1, \ldots, n\}$
so that two avoiding sets cannot occur. 

Following the classic papers \cite{TBA86, JLM03}, we may consider infinite left-products of stochastic matrices $\cdots Q(2) Q(1)$ such that: (i) the diagonal entries of all $Q(t)$ are positive (ii) all the positive entries of $Q(t)$ are bounded away from zero (ii) the graph obtained by unioning the edge sets of the graphs $G(Q(kB+1)), \ldots, G(Q((k+1)B))$ is strongly connected for all  $k$.  It is easy to see that, subject to these conditions, all the products $Q(kB+1) \cdots Q((k+1)B)$ will lie in some $P_{\gamma}$ with $\gamma>0$, and consequently the reasoning in the previous paragraph establishes the dynamic $x(t)=Q(t) x(t-1)$ converges to a multiple of the all-ones vector. 

\smallskip

\noindent {\bf Example 4:} It is possible to relax the connectivity conditions of the classic consensus results in several ways. For example, the results of \cite{morse-ieee} show the consensus property holds for a product $\cdots Q(2) Q(1)$ as long as each of the matrices $Q(t)$ has positive diagonal, the graphs $G(Q(t))$ all have a globally reachable node\footnote{Note that the globally reachable nodes may be different in different $G(Q(t))$.}, and the positive entries of $Q(t)$ are lower bounded by some $\gamma > 0$. A similar result is proved under more general
conditions in \cite{charon}. 

We may derive this from Theorem \ref{necsufftheorem} by arguing that these two conditions preclude the existence of two avoiding set sequences. Indeed,
letting $Q$ be a matrix satisfying the conditions in the previous paragraph, observe that the presence of self loops implies that $|N_Q(S)| \geq |S|$ for any subset $S$, since every node is its own neighbor. Moreover, the existence
of a globally reachable node implies that if $S,S'$ are disjoint and nonempty, we have that either $|N_Q(S)| > |S|$ or $|N_Q(S')| > |S'|$, since 
at least one of $S,S'$ does not contain the globally reachable node. This implies that if two avoiding set sequences $S_1, \ldots, S_l$ and $S_1', \ldots, S_l'$ exist, then $|S_i \cup S_i'|$ strictly increases whenever $i$ is increased to $i+1 ~{\rm mod}~ l$, which 
clearly can't be. 

\smallskip

\noindent {\bf Example 5:}  It was observed in \cite{morse-ieee} that a partial converse implication can also be shown to hold: if $\P$ is a compact set of matrices  and $\P$ is a consensus set, then each $P \in \P$ is has a globally reachable node. We next show how to deduce this quickly from Theorem \ref{necsufftheorem}. 

 Indeed, if $\P$ is a consensus set, then so is $\{P\}$ for any $P \in \P$. Consider the condensation of $G(P)$: this is directed acyclic graph obtained by collapsing together all the connected components of $G(P)$.  It is easy to see that the condensation of any graph has at least one sink node\footnote{In a directed graph, a sink node is a node with no out-going links.} and the condensation of $P$ must have exactly one sink node: if not, by choosing $S_1,S_1'$ corresponding to the connected components of two sink nodes, Theorem \ref{necsufftheorem} would imply that $\{P\}$ is not a consensus set.   Since the sink node is unique, it follows that all the nodes in the connected component corresponding to it are globally reachable.  

\smallskip

\noindent {\bf Example 6:} The previous examples are related to a classic theorem of Sarymsakov (\cite{sarym}, see also \cite{cao-sarym}). We will say that a stochastic matrix $P$ belongs to  the Sarymsakov
class ${\cal K}$ if it has the following property: for any two nonempty sets $A,B$ either $N_P(A)$ and $N_P(B)$ have an element in common, 
or $|N_P(A) \cup N_P(B) |> | A \cup B|$. Theorem \ref{necsufftheorem} then implies that any compact set of such matrices ${\cal K}$ is a consensus set, as the definition of ${\cal K}$ clearly precludes the existence of two 
avoiding set sequences: avoiding set sequences for matrices in ${\cal K}$ have the property that $|S_i \cup S_{i'}|$ increases whenever $i$ is increased to $i+1 ~{\rm mod}~ l$, which is impossible. 

\smallskip

\noindent {\bf Example 7:} A fairly sharp statement can be made for sets $\P$ of symmetric matrices. Namely, a compact set ${\cal P}$ of stochastic symmetric matrices is a consensus set if and only if $\{P\}$ is a consensus set for every $P \in \P$. Equivalently, a set $\P$ of symmetric matrices is a consensus set if and 
only if every matrix $P \in \P$ has a graph $G(P)$ which is connected and not bipartite. 
We next give a quick proof based on Theorem \ref{necsufftheorem}.

We argue that if the  avoiding set sequences  described  by Theorem \ref{necsufftheorem} exists, then at least one matrix $P \in \P$ has either a disconnected or bipartite graph $G(P)$. Indeed, observe that for symmetric $P$, we have that for any set $X$, $$ \sum_{i \in X, j \in \{1, \ldots, n\}} p_{ij} = \sum_{i \in \{1, \ldots, n\}, j \in X} p_{ij} = |X|.$$ Consequently,  picking any nonempty set $A$ and applying the above observation to both $X=A$ and $X=N_P(A)$, we see that we cannot have $|N_P(A)|<|A|$. Thus $|N_P(A)| \geq |A|$. Furthermore, again from the same observation,  $|N_P(A)| = |A|$ then $N_P(N_P(A)) = A$. 

Now if the sequences $S_i$ and $S_i'$ of Theorem \ref{necsufftheorem} exist, then the fact that for any $A$ we have $|N_P(A)| \geq |A|$ implies the cardinalities of all the $S_i$ are equal, 
as are the cardinalities of all the $S_i'$. Consequently, we have that in the graph $G(P(1)))$, $N(S_1)=S_2, N(S_2)=S_1$ and $N(S_1') = S_2', N(S_2')=S_1'$. 


If $|S_1 \cup S_2| < n$ or $|S_1' \cup S_2'| < n$ then $G(P)$ is not connected and we are finished. Else, we have that $S_1,S_2$ must be disjoint because if $a \in S_1 \cap S_2$ then 
we will not be able to have $S_1 \cap S_1' = \emptyset$ and $S_2 \cap S_2' = \emptyset$ since $a$ belongs to at least one of $S_1', S_2'$. The disjointness of $S_1,S_2$ along with the fact that $|S_1 \cup S_2| = n$ and $N(S_1)=S_2, N(S_2)=S_1$ imply $G(P(1))$ must be bipartite. 


\smallskip

\noindent {\bf Example 8:} We note that the claims of Example 7 cannot be generalized from symmetric to doubly stochastic matrices: it is not true that a compact 
set $\P$ of doubly stochastic matrices is a consensus set if and only if each matrix $P \in \P$ has the property that $\{P\}$ is a consensus set. Similarly,
it is not true that a compact set $\P$ of doubly stochastic matrices is a consensus set if and only if each matrix $P \in \P$ has a graph $G(P)$ which is strongly connected and not multipartite. 

A counterexample to both of these statements is the pair of matrices \[ A_1 = \left( \begin{array}{ccccc} 
0 & 1 & 0 & 0 & 0 \\
0 & 0 & 1/2 & 1/2 & 0 \\
0 & 0 & 1/2 & 0 & 1/2 \\
0 & 0 & 0 & 1/2 & 1/2 \\
1 & 0  & 0 & 0 & 0 
\end{array} \right), ~~~~~A_2 = \left( \begin{array}{ccccc} 0 & 0 & 0 & 0 & 1 \\ 
1 & 0 & 0 & 0 & 0 \\
0 & 1/2 & 1/2 & 0 & 0 \\
0 & 1/2 & 0 & 1/2 & 0 \\
0 & 0 & 1/2 & 1/2 & 0
\end{array} \right) \] Both $A_1,A_2$ are clearly doubly stochastic. By inspection, their graphs are strongly connected. Since  $(A_1)_{33} > 0, (A_2)_{33}>0$ neither graph is multipartite. Moreover, we claim that $\{A_1\}$ and $\{A_2\}$ are both consensus sets. This follows from the 
observation that in both $G(A_1)$ and $G(A_2)$ there is a path of length $4$ from any node to node $3$, so that two avoiding set sequences cannot occur. On the other hand, it is easy to verify that 
\[ S_1 = \{1\}, S_2 = \{2\} \]
\[ S_1' = \{5\}, S_2' = \{1\} \]
with the matrix sequence \[ A_1,A_2 \] are avoiding set sequences so that $\{A_1, A_2\}$ is
not a consensus set. 

\smallskip

\noindent {\bf Example 9:} The major new implication of the avoiding sets condition is the following fact, which is the main complexity-theoretic result of this section. 

\smallskip

\begin{proposition} Deciding whether a finite set of stochastic matrices is a consensus set is algorithmically decidable. 
\end{proposition} 

\smallskip

\begin{proof} By Theorem \ref{necsufftheorem}, this can be done by enumerating all pairs of set sequences in products of length at most $3^{n} - 2^{n+1} + 2$ from $\P$
and checking whether they avoid each other. The number of such set-sequences is clearly finite. 
\end{proof}

\subsection{Proof of Theorem \ref{necsufftheorem}} We will prove the necessicity and sufficiency of the avoiding sets condition in this subsection. First, we give a brief outline
of the proof.  We will find it considerably helpful to flip the direction of the products and 
consider infinite right-products rather than the infinite left-products in the definition of being a consensus set. In terms of infinite right-products, the property of being a consensus set 
can be seen to translate into a certain property about intersection of supports  - this is Lemma \ref{firstequivalence} below - which we will
show is equivalent to 
the ``avoiding sets'' property of the theorem. An old theorem of Paz \cite{paz} will provide an upper bound on the length $l$ of the sequence
until the two avoiding set sequences can be found. 

Our proof will make use of the following notion of ergodicity.

\smallskip

\begin{definition} For a stochastic matrix $P$ whose $i$'th row we will denote by $\p_i$, we define  \[ \delta(P) = \frac{1}{2} \max_{i,j=1,\ldots,n} ||\p_i-\p_j||_1.\]   The quantity $\delta(P)$ is called the Dobrushin ergodicity coefficient. It is not hard to see that for any two stochastic matrices $P_1, P_2$, we have \begin{equation} \label{smult} \delta (P_1 P_2) \leq \delta(P_1) \delta(P_2) \end{equation} See Chapter 6 of \cite{bremaud} for a proof. Since $0 \leq \delta(P) \leq 1$, we have in addition \begin{equation} \label{submult} \delta(P_1 P_2) \leq \min \left( \delta(P_1),  \delta(P_2) \right).\end{equation} \end{definition}

Our first lemma is an equivalence relation between a number of characterization of consensus sets. Its punchline is item \ref{inter}, which characterizes consensus sets in terms of the intersection of supports associated with infinite right products (rather than left products). 

\bigskip

\begin{lemma}  \label{firstequivalence} Let $\P$ be a compact set of stochastic matrices. The following are equivalent: \begin{enumerate} \item $\P$ is a consensus set. 
\item \label{equiv} For
every infinite sequence $P(1), P(2), \ldots$ from $\P$ the product
\[ \lim_{t \rightarrow \infty} P(t) P(t-1) \cdots P(1) \]
exists and has identical rows.
\item For every $\epsilon > 0$ there exists an integer $t(\epsilon)$ such that if $\Pi$ is the product of $t(\epsilon)$ matrices from 
$\P$, then $\delta(\Pi)<\epsilon$. 
\item \label{diff} For all $i,j$, and all infinite sequences $P(1), P(2), \ldots$ from ${\cal P}$, 
\[ \lim_{t \rightarrow \infty} (\e_i - \e_j) P(1) P(2) \cdots P(t) = 0.
\] Here  $\e_k$ stands for the $k$'th basis row vector.
\item \label{inter} There do not exist $i,j$ and an infinite sequence $P(1), P(2), \ldots$ from $\P$ such that 
\[ \e_i  P(1) P(2) \cdots P(t) \] and 
\[ \e_j  P(1) P(2) \cdots P(t) \] have disjoint 
supports for all $t$.

\item \label{lastitem} For all infinite sequences $P(1), P(2), \ldots$ from $\P$  and nonnegative vectors $\p_1, \p_2$ whose entries sum to $1$, 
\[ \lim_{t \rightarrow \infty} (\p_1 - \p_2) P(1) P(2) \cdots P(t) = 0. \] 
\end{enumerate}
\end{lemma}

\begin{proof} 

\bang{} {\em Equivalence of $(1)$ and $(2)$} {\ao is an obvious consequence of the definition of a consensus set.}

\smallskip

\bang{} {\em $(2)$ implies $(3)$: }%
Suppose item $(3)$ did not hold. That would mean there exists some $\epsilon > 0$ such that for every $i=1,2,\ldots$ there is a  product of at least $i$ matrices from $\P$ whose $\delta$ is at least $\epsilon$. 

Consider this infinite sequence of (finite) products, and in particular consider the right-most matrix in each product. Because $\P$ is compact, the right-most matrices have an accumulation point. Call this accumulation point $P_1$, and pick a subsequence of the products along which the right-most matrix converges to $P_1$. 

Next observe that, by the same argument, the second-right-most matrices in this new subsequence of products have an accumulation point. Call it $P_2$, and now further pick out a new subsequence (from the subsequence picked out in the previous paragraph) along which the second-right-most matrix converges to $P_2$. 

Proceeding in this way, we have an infinite sequence of stochastic matrices $P_1, P_2, P_3, \ldots$ with the following property: for any integer $k$, a product of arbitrarily small perturbations of the matrices $P_k, P_{k-1}, \ldots, P_1$ is the initial (right) part of a product whose $\delta$
is at least $\epsilon$.  More formally, given any integer $k$, there is a way to choose the matrices
$\Delta_1, \Delta_2, \Delta_3, \ldots, \Delta_k$ as close to zero as we like so that
\[ \delta \left( Q(t') \cdots Q(k+1) (P_k + \Delta_k) \cdots (P_3 + \Delta_3) (P_2 + \Delta_2) (P_1 + \Delta_1) \right) \geq \epsilon \] for some $t' \geq k+1$ and some stochastic matrices $Q(t'), \ldots,  Q(k+1)$.  But applying Eq. (\ref{submult}), this implies that 
\[ \epsilon   \leq \delta \left(  (P_k + \Delta_k) \cdots (P_3 + \Delta_3) (P_2 + \Delta_2) (P_1 + \Delta_1)  \right) \] Since $\delta$ is 
clearly a continuous function of the matrix elements, 
\[ \epsilon \leq \delta (P_k \cdots P_1) \] In turn, this implies that the infinite left-product 
\[ \cdots P_3 P_2 P_1   \] cannot approach a matrix with identical rows, which contradicts item $(2)$.

\smallskip

\bang{} {\em $(3)$ implies $(2)$:} Observe that item $(3)$ implies that for any sequence $P(1), P(2), \ldots $ from $\P$, the sequence $ \Pi_t $ of products $\Pi_t = P(t) \cdots P(1)$ is a Cauchy sequence. Indeed, for any $t_1,t_2 > t(\epsilon)$ the difference between a row of $\Pi_{t_1}$ and a row of $\Pi_{t_2}$ cannot be larger than $\epsilon$ in the $1$-norm because they are both convex combinations of the rows of the matrix $\Pi_{t(\epsilon)}$. Thus $\Pi_t$ has a limit. Now applying item $(3)$ once again we immediately get that this limit is a matrix with identical rows. 

\smallskip

\bang{} {\em Equivalence of $(3)$ and $(4)$:} That item $(3)$ implies $(4)$ is trivial. We next argue that if item $(3)$ is false, then item (4) is false as well. The proof is a variation of the argument we have used a few lines above to establish that item $(2)$ implied item $(3)$.   

Indeed, if item $(3)$ is false, then there exists $\epsilon > 0$ such
that for every $t$, there is a product of length at least $t$ of the matrices from ${\cal P}$ whose $\delta$ is at least $\epsilon$. Consider this infinite sequence of (finite) products, and in particular consider the left-most matrix in these products. Since $\P$ is compact, the left-most matrices have an accumulation point; call it $P_1$, and pick out a subsequence of these products along which the left-most matrix converges to $P_1$. 

Next, we consider the second-left-most matrix in this new subsequences of products. Since $\P$ is compact, the second-left-most matrices
have an accumulation point; call it $P_2$ and pick out (from the subsequence picked out in the previous paragraph) a new subsequence along which additionally the second-left-most 
matrix converges to $P_2$. 

Continuing in this way,
we obtain an infinite sequence of stochastic matrices $P_1, P_2, \ldots $ from ${\cal P}$ with the following property: given an integer $k$, there is a choice of 
the matrices $\Delta_1, \ldots, \Delta_k$ as close to zero as we like so that 
\[ (P_1 + \Delta_1) \cdots (P_k + \Delta_k) \] is the initial (left) part of a product whose $\delta$ is at least $\epsilon$. Appealing to 
Eq. (\ref{submult}) and continuity of $\delta$, we have that 
\[ \delta ( P_1 \cdots P_k) \geq \epsilon \] Adopting the notation $\Pi_k$ now for the right-product $P_1 \cdots P_k$, we have that for every $k$
there exist indices $j_1(k), j_2(k)$ such that 
\[ || \left( e_{j_1(k)} - e_{j_2(k)} \right) \Pi_k ||_1 \geq \epsilon \] Pick any pair $(j_1,j_2)$ which appears infinitely often in the set $\{ (j_1(k),j_2(k)) ~|~ k = 1, 2, \ldots \}$. 
We then have that 
\[ || \left( e_{j_1} - e_{j_2} \right) \Pi_k ||_1 \geq \epsilon \] for infinitely many $k$. This implies  item $(4)$ is false. 

\smallskip

\bang{} {\em Equivalence of $(4)$ and $(5)$:} It is easy to see that if item $(5)$ fails, then 
item $(4)$ fails. On the other hand, if item $(5)$ holds, it follows that for every pair of indices $i,j$ and sequence $P(1), P(2), \ldots$ from $\P$, there
is some integer $k$ (depending on $i,j$ and the sequence $P(1), P(2), \ldots$) such that rows $i$ and $j$ of $P(1) \cdots P(k)$ have positive entries in the same location. 
We claim this further implies that there is some finite integer $k'$ such that for {\bf all} sequences $P(1), P(2), \ldots $ of length $k'$, 
we have that {\bf every} pair of rows of the product $P(1) \cdots P(k')$ has a positive entry in the same location. 

This is true due by a variation of the argument we have deployed several times by now.  If no such $k'$ existed, then for every $t= 1,2, \ldots$ we could find a product of $t$ or more matrices from ${\cal P}$ with two rows whose supports do not intersect. We consider this infinite sequence of (finite) products and 
we pick a subsequence along which the left-most matrix converges to an accumulation point $P_1$; then, out of that subsequence, we pick a further subsequence out of which the second-left-most matrix converges to an accumulation point $P_2$, and so forth. We then obtain an infinite sequence $P_1, P_2, \ldots$ such that, for every $t$, there is a choice of matrices $\Delta_1, \ldots, \Delta_t$ as close to zero as we like with the product 
\[ (P_1 + \Delta_1) \cdots (P_t + \Delta_t) \] being the initial (left) part of a product with two rows that have disjoint supports. This means that 
\[ \delta \left(  (P_1 + \Delta_1) \cdots (P_t + \Delta_t) \right)  = 1 \] and consequently 
\[ \delta( P_1 \cdots P_t) = 1 \] and  the
product $P_1 \cdots P_t$ itself has two rows with disjoint supports. Picking now two rows which have disjoint
supports infinitely often in the sequence of products $$P_1,~~ P_1 P_2, ~~P_1 P_2 P_3,~~ \ldots$$ we obtain a contradiction to our initial statement, i.e., the existence of $k$ for every pair $i,j$ and matrix sequence. 

Now that we have established the existence of $k'$, let $\alpha$ be the infimum of $\delta(\Pi_{k'})$ as $\Pi_{k'}$ ranges over all possible products of length $k'$ from $\P$. We claim $\alpha < 1$. Indeed,
if $\alpha=1$, then by compactness of $\P$ we would be able to find a product of length $k'$ from $\P$ with two rows which have 
disjoint support. 

To summarize, we have shown that there exist an integer $k'$ and $\alpha \in [0,1)$ such that for every product $\Pi_{k'}$ from $\P$ of length $k'$, 
\[ \delta (\Pi_{k'}) \leq \alpha \] By Eq. (\ref{smult}), this means that for any positive integer $m$, 
\[ \delta ( \Pi_{m k'} )\leq \alpha^m \] and along with the monotonicity property of Eq. (\ref{submult}), this  implies item $(3)$ and consequently item $(4)$. 

\smallskip

\bang{} {\em Equivalence of $(4)$ and $(6)$:} Item $(4)$ is a special case of item $(6)$. To argue that item $(4)$ implies item $(6)$ observe that every vector whose entries sum to zero can be written as a linear combination of the vectors $\e_i - \e_j$. Thus for any $\p_1, \p_2$ satisfying the conditions of the lemma, there exist $a_{ij}$ such that \[ \p_1 - \p_2 = \sum_{i<j} a_{ij} (\e_i - \e_j),\] and so 
\[ \lim_{t \rightarrow \infty} (\p_1 - \p_2) P_{\tau(1)}  \cdots P_{\tau(t)} 
= \sum_{i<j} a_{ij} \lim_{t \rightarrow \infty} (\e_i - \e_j) P_{\tau(1)}  \cdots P_{\tau(t)}  = 0. \]
\end{proof}

\smallskip

\begin{remark} {\ao  We wish to remark on the following equivalence which was stated without proof in the introduction: 
${\cal P}$ is a consensus set if and only if for
every map $\tau: \Q \rightarrow \{1, \ldots, k\}$, the limit
\[ \lim_{t \rightarrow \infty} P(t) \cdots P(1)  \]
exists and has rank one. Indeed, if ${\cal P}$ is a consensus set then this statement immediately
follows from item $(2)$ of the preceeding lemma. Conversely, if $\lim_t P(t) \cdots P(1)  = uv^T$ for some vectors $u,v \in \R^n$, then because the set of stochastic matrices is a closed set which is also closed under multiplication, $uv^T$ must be stochastic as well; this forces $u = \1$ which implies item $(2)$ of the preceeding lemma.}  \end{remark}

\smallskip

With this lemma in place, we can prove Theorem \ref{necsufftheorem} without too much difficulty. The proof proceeds by relating the 
``two avoiding sets'' condition of the theorem to the ergodicity condition (namely item (5)) of the last lemma. 

\bigskip

\begin{proof}[Proof of Theorem \ref{necsufftheorem}] Suppose that such sequences of sets $S_1, \ldots, S_l$ , $S_1', \ldots, S_l'$ and matrices $P(1), \ldots, P(l)$ from $\P$ existed. Consider the infinite sequences of matrices whose first $l$ entries are $P(1), \ldots, P(l)$ and which repeats periodically thereafter.  Choose any $i \in S_1$ and $j \in S_1'$; we claim that 
the supports of $\e_i P(1) \cdots P(k) $ and $\e_j P(1) \cdots P(k)$ never intersect. 

Indeed, observe that the support of $\e_i P(1) \cdots P(k)$ is contained within $N^k(\{i\})$ which satisfies $N^k ( \{i\}) \subset S_{k +1 ~{\rm mod}~ l}$; similarly, the support of $\e_j P(1) \cdots P(k)$ is contained within $N^k (\{j\})$ which satisfies $N^k(\{j\}) \subset S_{k +1 ~{\rm mod}~ l}'$. Since by assumption $S_{k +1 ~{\rm mod}~ l} \cap S_{k +1~ {\rm mod}~ l}' = 
\emptyset$, we do indeed have that the supports of  $\e_i P(1) \cdots P(k)$ and $\e_j P(1) \cdots P(k)$ never
intersect. By item (5) of Lemma \ref{firstequivalence} we have that $\P$ is not a consensus set. 

This proves one direction. To prove the other, suppose now that no such sequences of sets $S_1, \ldots, S_l$, $S_1', \ldots, S_l'$ and matrices $P(1), \ldots, P(l)$ from $\P$ exist. Let $Q(1), Q(2), \ldots $ be an arbitrary sequence of stochastic matrices from ${\cal P}$; pick any $i$ and $j$ and define $U_k = N^k( \{i\})$ and $U_k' = N^k(\{j\})$. We claim that eventually $U_k$ and $U_{k}'$ constructed in this way intersect. Once this claim is proven, item (5) of Lemma \ref{firstequivalence} implies $\{P_1, \ldots, P_k\}$ is a consensus set, completing the proof.

To prove the claim, suppose instead that $U_k \cap U_{k}' = \emptyset$ for all $k$. Consider the collections of pairs of sets $(X,Y) \subset 2^{\{1, \ldots, n\}} \times 2^{\{1, \ldots, n\}}$. There are finitely many such pairs  so eventually
$(U_a,U_a') = (U_b,U_b')$ for some times $a,b$. By picking $S_1 = U_a, S_1' = U_a'$ and choosing \begin{eqnarray*} P(1) & = &  Q(a) \\ 
P(2) & = & Q(a+1) \\ 
\vdots & \vdots & \vdots \\ 
P(b-a) & = & Q( b-1) \\ 
\end{eqnarray*} and defining recursively \begin{eqnarray*} S_{k+1}  & = & N_{k} ( S_{k} ) \\
 S_{k+1}' & = & N_{k} (S_{k}') 
\end{eqnarray*} for $k=1, \ldots, b-a$, we obtain two ``avoiding set sequences.'' It remains to argue that the length $b-a$ may be upper bounded as in the theorem statement to obtain the desired contradiction. It is an old theorem of Paz (Theorem 4.7 in Chapter II.A of  \cite{paz}, proved in the course of showing the decidability of the ergodicity problem for inhomogenous Markov chains)  that 
the number of ordered partitions $(A,B,C)$ such that $A \cup B \cup C = \{1, \ldots, n\}$, $A \cap B = \emptyset, B \cap C = \emptyset, 
A \cap C = \emptyset$ and $A,B$ nonempty is $3^{n} - 2^{n+1}+1$. Consequently, if such two ``avoiding set cycles'' exist, by excising 
unnecessary repetitions we obtain they can be found
of length at most $3^{n} - 2^{n+1} + 1$. 
\end{proof}

\bigskip

Finally, we mention  another necessary and sufficient condition for consensus which follows out of the results of this section. Unlike the avoiding sets condition,
this condition is not immediately combinatorial and consequently less convenient to apply. On the other hand, verifying this condition algorithmically is a little easier than verifying the avoiding sets condition. 

\bigskip

\begin{corollary} A compact set $\P$ of $n \times n$ stochastic matrices is a consensus set if and only if every product of length $(1/2) (n-1) (3^n - 2^{n+1} + 1)$ of matrices from $\P$ has a positive column.
\end{corollary} 

\bigskip

Note that because the limit of every infinite left-product from a consensus set is a stochastic matrix with identical rows, we have that every product from a consensus set eventually has
a positive column.  The novelty here comes from the explicit upper bound on the length of the product. 

\bigskip

\begin{proof} One direction is easy - if every product of length $(1/2) (n-1)(3^n - 2^{n+1} + 1)$ from $\P$ has a positive column, then it is well-known that $\P$ is a consensus set: see Example 2 in the
previous Section \ref{examples}. 

Conversely, suppose $\P$ is a consensus set.  Item $(5)$ of Lemma \ref{firstequivalence} implies that for every sequence  $P(1), P(2), \ldots$ from $\P$ and any pair $i,j \in \{1, \ldots, n\}$,
the supports of the vectors $\e_i P(1) P(2) \cdots P(t)$ and $\e_j P(1) \cdots P(t)$ eventually intersect for some $t$ large enough. Recall that we showed the course of the proof of item $(5)$ that this
in turn
implies that  there is a uniform upper bound 
$k$ such that {\bf every} possible $\e_i P(1) \cdots P(t)$ and $\e_j P(1) \cdots P(t)$ have supports which intersect whenever $t \geq k$. 

We claim this now further implies that the product of every $(n-1)k$ matrices from $\P$ has a positive column. We argue this next by inductively proving the more general statement that, for every $m=1, \ldots,n-1$, the product of 
every $mk$ matrices from $\P$ has a column which has at least $m+1$ nonzero entries.

 Indeed, we have already established the base case of $m=1$. Now suppose we have proved this for some $m$;
that is to say, suppose we have shown that for every product $P(1) \cdots P(mk)$ from $\P$, there is an index $c \in \{1, \ldots, n\}$ (depending on the product) such that the $c$'th column of $P(1) \cdots P(mk)$ 
has nonzero entries in some set $N \subset \{1,\ldots,n\}$ with cardinality $|N|=m+1$. Pick any index $i$ which is not in $N$ and observe that since the supports of  $\e_c P(mk+1) \cdots P((m+1)k)$ and $\e_i P(mk+1) \cdots P((m+1)k)$ intersect, it follows that there is a column of $P(1) \cdots P(mk+1)$ 
which has nonzero values in the entries of $N \cup \{i\}$. This proves the inductive step and thus proves the claim that every product of $(n-1)k$ matrices from $\P$ has a positive column.

Finally, it is a result of Paz (see Theorem 4.7 in Section II.A of \cite{paz}) that if there exists some $l$ such that every product of $l$ matrices from a {\em finite} set $\P$ of $n \times n$ stochastic matrices
is scrambling, then one may take $l = (1/2) (3^{n} - 2^{n+1} + 1)$. Since the number of matrices that appear in any single product of length at most $k$ is finite, Paz's result 
implies the current corollary. 
\end{proof}
  
  
  \section{Complexity of the consensus problem\label{hard}} In this section, we show that unless $P=NP$ there does not exist a polynomial-time algorithm for the consensus problem even in the case of two matrices. We prove a similar result for three undirected matrices. We begin by stating our main results.

\bigskip

\begin{theorem} \label{twotheorem} Unless $P=NP$, the problem of determining if a given set of \emph{two} stochastic matrices is a consensus set cannot be decided in polynomial-time.  \end{theorem}

\bigskip

\begin{theorem} \label{threetheorem} 
Unless $P=NP$, the problem of determining if a given set of \emph{three undirected} stochastic matrices is a consensus set cannot be decided in polynomial-time.
 \end{theorem}

  \bigskip
  
We remark that that the assumptions of Theorem \ref{threetheorem} cannot be relaxed from requiring each of the matrices to be undirected to requiring the matrices to be symmetric. Indeed, as Example 7 in Section \ref{examples} shows, it is possible to check 
whether a finite collection of symmetric stochastic matrices is a consensus set in polynomial time. 

\bigskip

The remainder of this section will be devoted to the proof of these theorems, which we now briefly outline. 

Our strategy will be to show both consensus and the non-satisfiability of a $3$-SAT problem are equivalent to the non-existence of a path between two nodes in every possible sequence of two particularly chosen graphs $G_1, G_2$ (our notion of a path in a graph sequence is the straightforward one; the interested reader may refer to Definition \ref{pathdef} below for a precise statement).  One side of this equivalence is easy: the equivalence of consensus and the non-existence of a path is easily obtained by taking scalings of the adjacency matrices of the graphs and applying the results of the previous section on supports which eventually overlap. On the other hand, to show that the unsatisfiability of $3$-SAT and the non-existence of a path are equivalent we use a variation of a previous
construction from \cite{TB97}. Finally, to prove the result for undirected matrices we will construct a gadget for simulating directed paths on sequences of two graphs with undirected paths on sequences of three graphs. 

We next proceed to the details of the proofs. Our first step is to define the $3$-SAT problem. 

\bigskip

\begin{definition} A clause is defined to be the logical OR of Boolean variables or their negations. For example, $x_1 \vee x_5 \vee \overline{x}_7 \vee \overline{x}_8$ is a clause. The variables $x_i$ and their negations $\overline{x}_i$ are called 
litorals. A $3$-SAT formula is the logical AND of clauses which contain $3$ litorals. For example, $(x_1 \vee x_2 \vee x_3) \wedge (\overline{x}_1 \vee 
\overline{x}_2 \vee \overline{x}_3)$ is a a $3$-SAT formula. A $3$-SAT formula on Boolean variables $x_1, \ldots, x_n$ is called
satisfiable if is possible to assign a $0$ or a $1$ to each $x_i$ so that the formula evaluates to $1$. For example, the formula  $(x_1 \vee x_2 \vee x_3) \wedge (\overline{x}_1 \vee 
\overline{x}_2 \vee \overline{x}_3)$ is satisfiable and $x_1 = 0, x_2 = 0, x_3=1$ is one assignment which satisfies it. The $3$-SAT problem is the problem of 
determining if a given formula is satisfiable. 
\end{definition} 

\bigskip We next describe the graphs which will provide the ``bridge'' between $3$-SAT and consensus problems.  

\bigskip

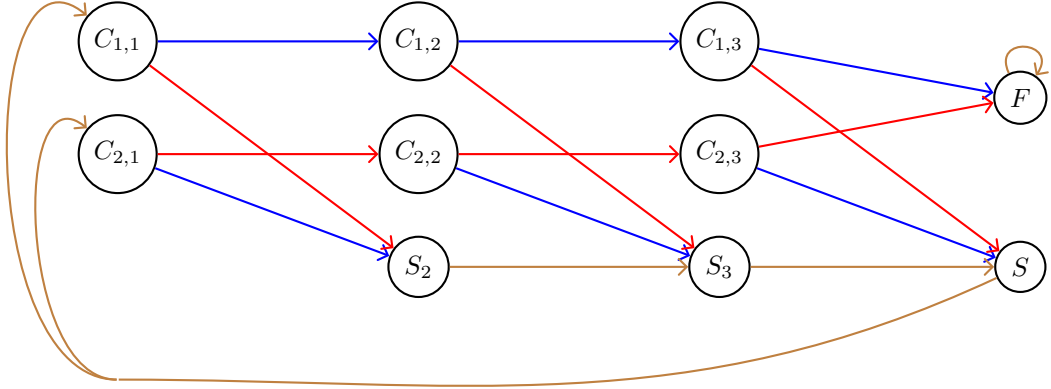
\begin{figure} \begin{center}
\begin{tikzpicture}[->, thick]
\SetVertexNormal[LineColor=black]
\SetVertexMath

\node (T1) at (0,3) [circle, draw] {$C_{1,1}$};
\node (T2) at (4,3) [circle, draw] {$C_{1,2}$};
\node (T3) at (8,3) [circle, draw] {$C_{1,3}$};
\node (S) at (12,2.25) [circle, draw] {$F$};
\node (M1) at (0,1.5) [circle, draw] {$C_{2,1}$};
\node (M2) at (4,1.5) [circle, draw] {$C_{2,2}$};
\node (M3) at (8,1.5) [circle, draw] {$C_{2,3}$};
\node (B2) at (4,0) [circle, draw] {$S_2$};
\node (B3) at (8,0) [circle, draw] {$S_3$};
\node (F) at (12,0) [circle, draw] {$S$};


\node (Q1) at (0,-1) [circle, draw, minimum size=0pt, inner sep=0pt] {};
\node (Q2) at (0,-1.5) [circle, draw, minimum size=0pt, inner sep=0pt] {};

\path[->,color=blue, >=angle 90]
(T1) edge (T2)
(T2) edge (T3)
(T3) edge (S)

(M1) edge (B2)
(M2) edge (B3)
(M3) edge (F);



\path[->,color=red, >=angle 90]

(T1) edge (B2)
(T2) edge (B3)
(T3) edge (F)

(M1) edge (M2)
(M2) edge (M3)
(M3) edge (S);






\path[->,color=brown, >=angle 90]
(B2) edge  (B3)
(B3) edge  (F)
(F) edge[out=205, in=0,-] (Q2)
(Q2) edge[out=180, in=140] (T1)
(Q2) edge[out=180, in=140] (M1)
(S) edge [out=115, in=55, looseness=4] (S);


\end{tikzpicture} \caption{The graphs corresponding to the formula $( x_1 \vee x_2 \vee x_3) \wedge (\overline{x_1} \vee \overline{x_2} \vee \overline{x_3})$. Edges which are present in both $G_0(f), G_1(f)$ are shown in brown; edges which are only in $G_0(f)$ are shown in blue and edges which are only in $G_1(f)$ are in red.  \label{twoconsensus}}  \end{center} \end{figure}

\begin{definition} Given a $3$-SAT formula $f$, we will define two graphs $G_0(f), G_1(f)$ as follows (the reader may wish to refer to Figure \ref{twoconsensus} for an example).

We will create a node for each clause/variable pair; the node corresponding to clause $i$ and variable $j$ will
be called $C_{i,j}$. There will be a node
$F$ (intuitively thought of as a ``failure node'') and nodes $S_2, \ldots, S_{n+1}$ (intuitively thought of as ``intermediate success nodes.''). 

For convenience, we will also adopt the following notation. For any $i$, $C_{i,n+1}$ will refer to the node $F$. Moreover, $S_{n+1}$ will sometimes be labeled simply as $S$.

If setting variable $j$ to $0$ satisfies clause $i$, then we will put a directed link from $C_{i,j}$ to $S_{j+1}$ in $G_0(f)$; else, we will put a directed link from $C_{i,j}$ to $C_{i,j+1}$ in $G_0(f)$ (note that when $j=n$, this link leads to $F$). 

If setting variable $j$ to $1$ satisfies clause $i$, then we will put a directed link from $C_{i,j}$ to $S_{j+1}$ in $G_1(f)$; else, we will put a directed link from $C_{i,j}$ to $C_{i,j+1}$ in $G_1(f)$ (note that when $j=n$ this link leads to $F$). 

Finally, both $G_1,G_2$ will have self-loops at $F$, links from each $S_{i}$ to $S_{i+1}$ and links from $S=S_{n+1}$ to all 
$C_{i,1}$. \end{definition}

\bigskip

We will require the notion of a path in a time-varying graph sequence. 

\bigskip

\begin{definition} \label{pathdef} Consider a sequence of {\ao directed} graphs $G_1, G_2,  \ldots, G_t$ on the same vertex set $V$. We will say that $e_{1}, \ldots, e_{t}$ is a path from node $u$ to
node $v$ in this sequence if:  (i) $e_{i}$ is an edge in $G_{i}$ {\ao(ii) the source of $e_1$ is $u$ and the destination of $e_{t}$ is $v$ (iii) the source of $e_{i}$ is the same as the destination of $e_{i-1}$ for $i=2,\ldots,t$.  }
\end{definition} 

\bigskip With these definitions in place, we can proceed to our first lemma below. Our first step is to relate satisfiability of the formula $f$ to the existence of
a path in a switching of the graphs $G_0(f),G_1(f)$, and this is what next lemma accomplishes. 

\bigskip

\begin{lemma} \label{pathhardness} A $3$-SAT formula $f$ is satisfiable if and only if there exists a sequence of
$G_0(f), G_1(f)$ of length at least $n+1$ without a path from $S$ to $F$. \end{lemma}

\begin{proof} Suppose that the $3$-SAT problem is satisfiable, i.e., there exists an assignment of $\{0,1\}$ to the variables satisfying every clause. We need to find a sequence of the graphs $G_0(f), G_1(f)$ of length $n+1$ without a path from $S$ to $F$. Because the links going out of node $S$ are the same in both $G_0(f)$ and $G_1(f)$ (they lead to all of the nodes $C_{i,1}$), equivalently we need to find a sequence of the graphs $G_0(f), G_1(f)$ of length $n$ such that there is no path from any $C_{i,1}$ to $F$. 

We will construct such a sequence as follows: we will pick a satisfying assignment and we will set the $i$'th graph to be $G_0(f)$ if the $i$'th variable in the satisfying assignment is $0$ and $G_1(f)$ if the $i$'th variable in the satisfying assignment is $1$. 

We argue there is no path from any $C_{i,1}$ to $F$ in this sequence. Indeed, the links from any $C_{i,j}$ either lead ``to the right'' to $C_{i,j+1}$ or ``down'' to $S_{j+1}$. Because there is no path of length  $n$ or less from any $S_{j}$ to $F$, it suffices to argue that at some time step, only the ``down'' edge to some $S_{j+1}$ 
will be available. But clause $i$ is satisfied by at least one variable in the satisfying assignment, say by variable $k$; and so by construction the only edge from $C_{i,k}$ at time $k$ will lead to $S_{k+1}$. Thus there is no path from $C_{i,1}$ to $C_{i,k+1}$ of length $k$, which means that every path beginning at $C_{i,1}$ of length $k$ arrives at $S_{k+1}$. This proves the ``only if'' part. 

To prove the ``if'' part, assume that there is a sequence of length $n+1$ with no path from $S$ to $F$; as we remarked above, this is equivalent to the assumption that there a sequence of length $n$ with no path from
all the $C_{i,1}$ to $F$. Let $x_i$ be $0$ if the $i$'th element of the sequence is $G_0(f)$ and $1$ if the $i$'th element of the sequence is $G_1(f)$. We argue that this is a satisfying assignment. 

Indeed, consider some clause $i$. We know that there is no path of length $n$ from $C_{i,1}$ to $F$. Note that if at every time step $k$ there was a link from $C_{i,k}$ to $C_{i,k+1}$ then a path from $C_{i,1}$ to $F$ would exist. Thus  there is some time step $k$ at which there is only a link from $C_{i,k}$ to $S_{k+1}$.  By construction of the graphs $G_0(f), G_1(f)$,  this means that our assignment satisfies clause $i$ via variable $x_k$.  \end{proof}

\bigskip

%
%

\begin{definition} \label{amatrix} Define $d_{0}(f,i)$ to be the out-degree of node $i$ in $G_0(f)$; $d_1(f,i)$ is defined similarly. We will define the matrix $A_0(f)$ as $[A_0(f)]_{ij} = 1/d_0(f,i)$ if $(i,j)$ is an edge in $G_0(f)$ and $0$ otherwise;  the matrix $A_1(f)$ is defined similarly. Note that $A_0(f)$ and $A_1(f)$ are nonnegative stochastic matrices.
\end{definition}

\bigskip

\begin{lemma}\label{pathlemma} $\{A_0(f), A_1(f)\}$ is a consensus set if and only if every sequence of graphs $G_0(f), G_1(f)$ of length $n+1$ has a path from $S$ to $F$. \end{lemma}

\begin{proof} By item $(5)$ of Lemma \ref{firstequivalence},  $\{A_0(f), A_1(f)\}$ is a consensus set if and only if  for every sequence of $G_0(f), G_1(f)$  and every pair of indices $i,j$, there is a third index $k$ reachable from both $i$ and $j$ through a path in that graph sequence.

This, however, is equivalent to the assertion that every graph sequence of length $n+1$ has a path from $S$ to $F$. Indeed, if every graph sequence of length $n+1$ had such a path, then from every node $k$ (and for every infinite graph sequence) there is a path from  $k$ to $F$ of length $2(n+1)$ (by going from $k$ to one of $\{S,F\}$ in at most $n+1$ steps and then either taking the selp-loop at $F$ or going from $S$ to $F$ in $n+1$ more). 

Conversely, if there is a sequence of length $n+1$ with no path from $S$ to $F$, then by concatenating it with itself we would obtain an infinite graph sequence with the property that the set of nodes reachable from $F$ (which is just $F$) and the set of nodes reachable from $S$ never intersect.  \end{proof}

\bigskip

%

\begin{remark}  \label{pathn1} The above proof also shows that an infinite sequence of $\{G_0(f), G_1(f)\}$ without a path from $S$ to $F$ exists if and only if such a sequence of length
$n+1$ exists. Consequently, $\{A_0(f), A_1(f) \}$ is a consensus set if and only if every infinite sequence of graphs $G_0(f), G_1(f)$ has a path from
$S$ to $F$. We remark on this here  because we will use this property later on in the proof of the subsequent Theorem \ref{threetheorem}. 
\end{remark}

\bigskip

\par{\it Proof of Theorem \ref{twotheorem}}:  By putting together Lemma \ref{pathhardness} and Lemma \ref{pathlemma} we observe that, starting from an instance of $3$-SAT with variables $x_1, \ldots, x_n$ and $m$ clauses, we can construct in polynomial time nonnegative, stochastic matrices $P_1,P_2 \in \R^{(m+1)n{\ao + 1} ~\times~ (m+1)n+{\ao 1}}$ such that the $3$-SAT instance is unsatisfiable if and only if $\{P_1, P_2\}$ is a consensus set. The theorem immediately follows.

{\vbox{\hrule height0.6pt\hbox{%
   \vrule height1.3ex width0.6pt\hskip0.8ex
   \vrule width0.6pt}\hrule height0.6pt
  }}
  
  \bigskip
  
  \begin{remark} Observe that for any formula $f$, an infinite left-product of of the matrices $\{A_0(f), A_1(f)\}$ either does not converge to a rank-$1$ matrix or converges to the matrix $\e_F {\bf 1}^T$. It may be argued that this represents a ``pathological'' case of convergence to consensus as in many applications it is desired that the final matrix equal $v \1^T$ where $v$ is a strictly positive vector; intuitively, one wants every node to "contribute'' to the final result.

  We give a brief informal sketch here of a construction that remedies this issue. We will create two new matrices $A_0(f), A_1(f)$ with the property that they form a consensus set if and only if old $A_0(f), A_1(f)$ formed a consensus set; however, any convergent infinite left-product of the new matrices must converge to a strictly positive matrix (which, consequently, can only be written as $v {\bf 1}^T$ with strictly positive $v$). These matrices are constructed exactly as before in Definition \ref{amatrix} from, however, new graphs $G_0(f), G_1(f)$ defined as follows. The new $G_0(f)$ is obtained by taking two copies of the old $G_0(f)$, adding links going from each of the $F$ nodes to all the nodes in its copy, and adding links going both ways between the two $F$ nodes. The new $G_1(f)$ is produced in 
  the same way. We recommend the reader consult Figure \ref{doubledfig} for an illustration.

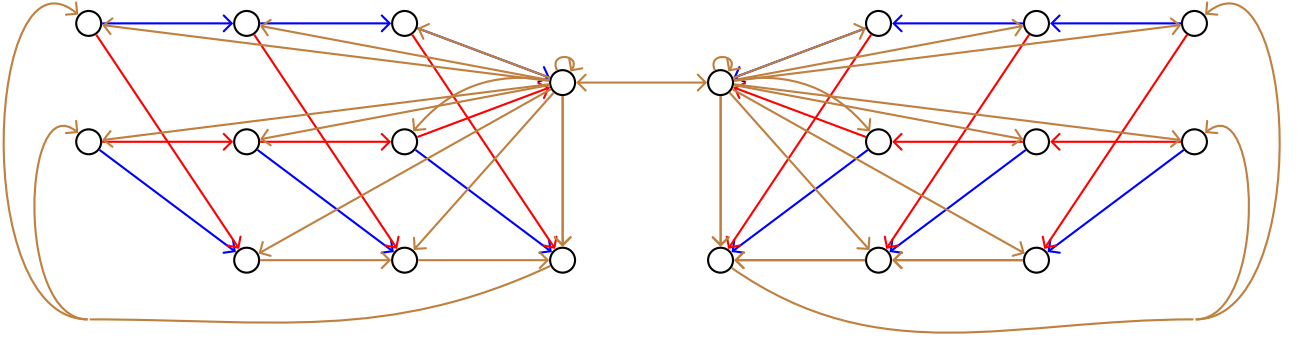
\begin{figure}
\begin{tikzpicture}[->, thick,scale=1.05]
\SetVertexNormal[LineColor=black]
\SetVertexMath

\node (T1) at (0,3) [circle, draw] {};
\node (T2) at (2,3) [circle, draw] {};
\node (T3) at (4,3) [circle, draw] {};
\node (S) at (6,2.25) [circle, draw] {};
\node (M1) at (0,1.5) [circle, draw] {};
\node (M2) at (2,1.5) [circle, draw] {};
\node (M3) at (4,1.5) [circle, draw] {};
\node (B2) at (2,0) [circle, draw] {};
\node (B3) at (4,0) [circle, draw] {};
\node (F) at (6,0) [circle, draw] {};

\node (T1right) at (14,3) [circle, draw] {};
\node (T2right) at (12,3) [circle, draw] {};
\node (T3right) at (10,3) [circle, draw] {};
\node (Sright) at (8,2.25) [circle, draw] {};
\node (M1right) at (14,1.5) [circle, draw] {};
\node (M2right) at (12,1.5) [circle, draw] {};
\node (M3right) at (10,1.5) [circle, draw] {};
\node (B2right) at (12,0) [circle, draw] {};
\node (B3right) at (10,0) [circle, draw] {};
\node (Fright) at (8,0) [circle, draw] {};


\node (Q1) at (0,-1) [circle, draw, minimum size=0pt, inner sep=0pt] {};
\node (Q2) at (0,-0.75) [circle, draw, minimum size=0pt, inner sep=0pt] {};
\node (Q2right) at (14,-0.75) [circle, draw, minimum size=0pt, inner sep=0pt] {};

\path[->,color=blue, >=angle 90]
(T1) edge (T2)
(T2) edge (T3)
(T3) edge (S)

(M1) edge (B2)
(M2) edge (B3)
(M3) edge (F)

(T1right) edge (T2right)
(T2right) edge (T3right)
(T3right) edge (Sright)

(M1right) edge (B2right)
(M2right) edge (B3right)
(M3right) edge (Fright);

\path[->,color=red, >=angle 90]

(T1) edge (B2)
(T2) edge (B3)
(T3) edge (F)

(M1) edge (M2)
(M2) edge (M3)
(M3) edge (S)

(T1right) edge (B2right)
(T2right) edge (B3right)
(T3right) edge (Fright)

(M1right) edge (M2right)
(M2right) edge (M3right)
(M3right) edge (Sright);

\path[->,color=brown, >=angle 90]

(S) edge (F)
(S) edge (T1)
(S) edge (T2)
(S) edge (T3)
(S) edge (B2)
(S) edge (B3)
(S) edge (M1)
(S) edge (M2)
(S) edge [bend right] (M3)
(S) edge (F)
(B2) edge  (B3)
(B3) edge  (F)

(Sright) edge (Fright)
(Sright) edge (T1right)
(Sright) edge (T2right)
(Sright) edge (T3right)
(Sright) edge (B2right)
(Sright) edge (B3right)
(Sright) edge (M1right)
(Sright) edge (M2right)
(Sright) edge [bend left] (M3right)
(Sright) edge (Fright)
(B2right) edge  (B3right)
(B3right) edge  (Fright)

(F) edge[out=205, in=0,-] (Q2)
(Q2) edge[out=180, in=140] (T1)
(Q2) edge[out=180, in=140] (M1)
(S) edge [out=115, in=55, looseness=4] (S)

(B2right) edge  (B3right)
(B3right) edge  (Fright)
(Fright) edge[out=-35, in=180,-] (Q2right)
(Q2right) edge[out=0, in=40] (T1right)
(Q2right) edge[out=0, in=40] (M1right)
(Sright) edge [<->] (S)
(Sright) edge [out=115, in=55, looseness=4] (Sright);


\end{tikzpicture} \caption{The graphs corresponding to the formula $( x_1 \vee x_2 \vee x_3) \wedge (\overline{x_1} \vee \overline{x_2} \vee \overline{x_3})$. Edges which are present in both $G_0(f), G_1(f)$ are shown in brown; edges which are only in $G_0(f)$ are shown in blue and edges which are only in $G_1(f)$ are in red.  \label{doubledfig}}  \end{figure}
  
  Observe that if the old $\{A_0(f), A_1(f)\}$ was not a consensus set, the new matrices do not form a consensus set either. Indeed, as we have argued the old set of matrices was a consensus set if and only if there existed a sequence of the old $G_0(f), G_1(f)$ such that starting from $S$ one never reached $F$; this implies that in the same sequence of the new $G_0(f), G_1(f)$, paths starting from the two $S$ nodes never reach their corresponding $F$ nodes and so never intersect. This implies the new matrices $A_0(f), A_1(f)$ do not form a consensus set. 
  
  Conversely, if the old $\{A_0(f), A_1(f)\}$ was a consensus set, the new matrices form a consensus set as well. By Lemma \ref{pathlemma}, we know that that the old matrices formed a consensus set if and only if in any sequence of graphs there is a path of length $2(n+1)$ from each node to node $F$ (by going first to one of $\{S,F\}$ in at most $n+1$ steps and then either taking the self-loop at $F$ $n+1$ times or going from $S$ to $F$). It immediately follows that in any sequence of the new graphs, there is a path from any node to any other node of length $2(n+1)+2$. By item (5) of Lemma \ref{firstequivalence} this implies that 
$\{A_0(f), A_1(f)\}$ is a consensus set; moreover, it implies that any product of the matrices $A_0(f), A_1(f)$  of length $2(n+1)+2$ has
entries bounded below by $(1/(2(m+1)(n+1)))^{2(n+1)+2}$, which means that any limiting matrix of an infinite left-product must be positive.  
   \end{remark}

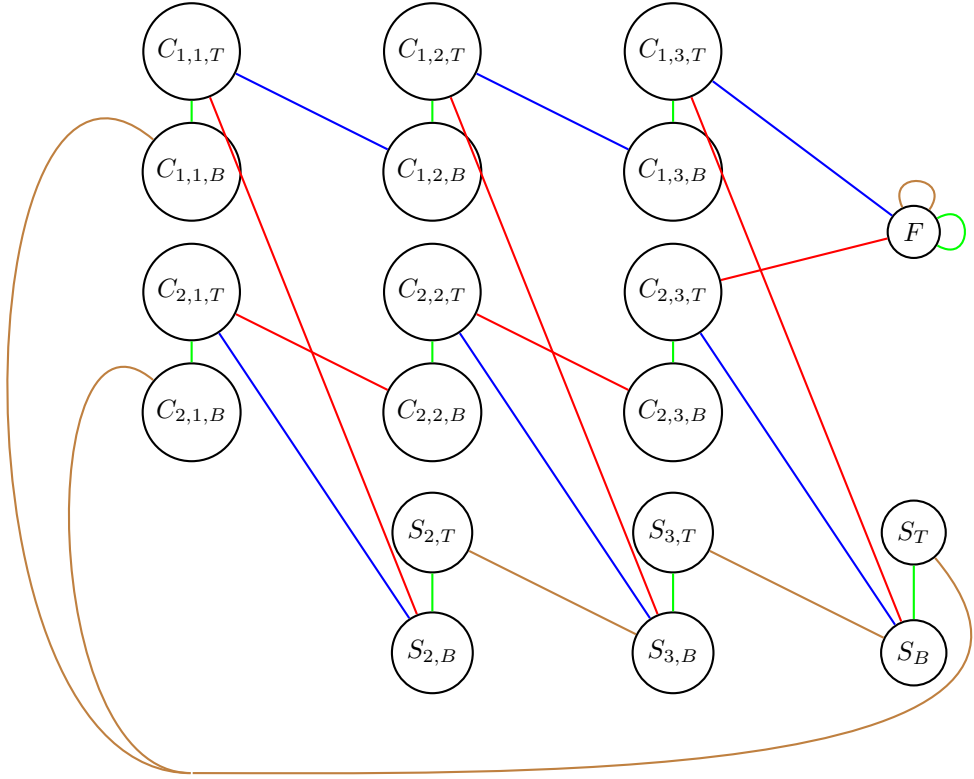
\begin{figure} \begin{center}
\begin{tikzpicture}[thick,scale=0.8]
\SetVertexNormal[LineColor=black]
\SetVertexMath

\node (T1) at (0,8) [circle, draw] {$C_{1,1, B}$};
\node (T12) at (0,10) [circle, draw] {$C_{1,1, T}$};

\node (T2) at (4,8) [circle, draw] {$C_{1,2,B}$};
\node (T22) at (4,10) [circle, draw] {$C_{1,2,T}$};

\node (T3) at (8,8) [circle, draw] {$C_{1,3,B}$};
\node (T32) at (8,10) [circle, draw] {$C_{1,3,T}$};

\node (S) at (12,7) [circle, draw] {$F$};

\node (M1) at (0,4) [circle, draw] {$C_{2,1,B}$};
\node (M12) at (0,6) [circle, draw] {$C_{2,1,T}$};

\node (M2) at (4,4) [circle, draw] {$C_{2,2,B}$};
\node (M22) at (4,6) [circle, draw] {$C_{2,2,T}$};

\node (M3) at (8,4) [circle, draw] {$C_{2,3,B}$};
\node (M32) at (8,6) [circle, draw] {$C_{2,3,T}$};

\node (B2) at (4,0) [circle, draw] {$S_{2,B}$};
\node (B22) at (4,2) [circle, draw] {$S_{2,T}$};

\node (B3) at (8,0) [circle, draw] {$S_{3,B}$};
\node (B32) at (8,2) [circle, draw] {$S_{3,T}$};

\node (Ftop) at (12,2) [circle, draw] {$S_T$};
\node (Fbottom) at (12,0) [circle, draw] {$S_B$};


\node (Q1) at (0,-1) [circle, draw, minimum size=0pt, inner sep=0pt] {};
\node (Q2) at (0,-2) [circle, draw, minimum size=0pt, inner sep=0pt] {};

\path[color=blue, >=angle 90]
(T12) edge (T2)
(T22) edge (T3)
(T32) edge (S)

(M12) edge (B2)
(M22) edge (B3)
(M32) edge (Fbottom);

\path[color=red, >=angle 90]

(T12) edge  (B2)
(T22) edge  (B3)
(T32) edge  (Fbottom)

(M12) edge (M2)
(M22) edge (M3)
(M32) edge (S);



\path[color=green, >=angle 90]

(T1) edge (T12)
(T2) edge (T22)
(T3) edge (T32)
(M1) edge (M12)
(M2) edge (M22)
(M3) edge (M32)
(B2) edge (B22)
(B3) edge (B32)


(Ftop) edge (Fbottom)
(S) edge [out=30, in=-30, looseness=4] (S);

\path[color=brown, >=angle 90]
(S) edge [out=115, in=55, looseness=4] (S)
(Ftop) edge[out=310, in=0,-] (Q2)
(Q2) edge[out=180, in=140] (T1)
(Q2) edge[out=180, in=140] (M1)
(B22) edge  (B3)
(B32) edge  (Fbottom);


\end{tikzpicture} \vspace{-10pt} \caption{The graphs corresponding to the formula $( x_1 \vee x_2 \vee x_3) \wedge (\overline{x_1} \vee \overline{x_2} \vee \overline{x_3})$. The common edges of $s_0(f), s_1(f)$ are in brown; the edges which are only in $s_0(f)$ are in blue, while the edges which are only in $s_1(f)$ are in red. The graph $s_2(f)$ is given in green. \label{threeconsensus} } \end{center} \end{figure}

\bigskip

\par{\it Proof of Theorem \ref{threetheorem}}: Given a $3$-SAT formula $f$, we will construct three  
undirected stochastic matrices $B_0(f), B_1(f), B_2(f)$ whose dimensions are polynomial in $m,n$ and we will prove that $\{ B_0(f), B_1(f), B_2(f) \}$ is a 
consensus set if and only if $\{A_0(f), A_1(f) \}$ is a consensus set. Together with Theorem 
\ref{twotheorem} this proves the current theorem. 

We will construct these matrices as in Definition \ref{amatrix} from undirected graphs $g_0(f), g_1(f), g_2(f)$
defined next. Note that the undirectedness of the graphs will imply the undirectedness of the corresponding 
matrices. 

\smallskip

\noindent {\bf The construction:} We will first construct undirected graphs $s_0(f), s_1(f), s_2(f)$ by taking every node of $G_0(f),G_1(f)$ (recall that these graphs share a common
set of nodes) except $F$ and splitting them into two nodes - a ``top'' node and ``bottom'' node.  We will use the subscripts $T,B$ to refer to the new nodes, e.g., node $a \neq F$ splits into a top node $a_T$ and a bottom
node $a_B$. For convenience of notation, both $F_T$ and $F_B$ will refer to the (single) node $F$. We will put edges in the new (undirected) graphs $s_0(f), s_1(f)$ based on old edges in the (directed) graphs $G_0(f), G_1(f)$ as follows: if $(k,l)$ was an edge in $G_i(f)$ then we will put the 
edge $(k_T, l_B)$ in $s_i(f)$. The graph $s_3(f)$ will contain all edges from the top node of a split pair to the bottom node of the same split pair (which we will understand to include a self-loop at node $F$). We refer the reader to Figure \ref{threeconsensus} for a concrete example of these graphs. 

We finally construct $g_0(f), g_1(f)$ as follows. We take $s_0(f),s_1(f)$ and add edges from every bottom node to $F$; then we take two copies 
of the resulting graph and connect their $F$ nodes. The graph $g_3(f)$ is constructed in the same way, except that the new edges in each copy go from every top node to $F$ instead. 

We remark that the matrices $B-(f), B_1(f), B_2(f)$ which are constructed from Definition \ref{amatrix} with graphs $g_0(f), g_1(f), g_2(f)$ are stochastic matrices because every node in the latter graphs is incident on at least one edge.  
\smallskip

\noindent {\bf Properties of the construction:} Our construction has the following property: if there is a path from node $a$ to node $b$ of length $k$ in some  sequence of $G_0(f), G_1(f)$ then in each of the two copies there is a path from $a_T$ to $b_T$  of length $2k$ in the graph sequence obtained by replacing $G_0(f)$ by $g_0(f)$, $G_1(f)$ by $g_1(f)$, and inserting $g_2(f)$ into every even 
time step. For example, if there is a path from $a$ to $b$  of length $4$ in the sequence 
\[  G_0(f) G_0(f) G_1(f) G_0(f) \cdots \] then there is a path from $a_T$ to $b_T$ of length $8$ in the sequence 
\[ g_0(f) g_2(f)  g_0(f) g_2(f) g_1(f) g_2(f) g_0(f) g_2(f) \cdots \] Thus paths in graph sequences of $G_0(f), G_1(f)$ can be ``ported'' into
sequences of the graphs $g_0(f), g_1(f), g_2(f)$ wherein $g_2(f)$ appears at, and only at, the even  time steps.

Next we remark that a partial converse of this statement is also true:  if there is a path from $a_T$, which is not one of the two $F$ nodes, to the $F$ node in its copy of length $2k$ in a sequence of $g_0(f), g_1(f), g_2(f)$ wherein $g_2(f)$ appears at, and only at, all the even time steps, then there is a path from 
$a$ to $F$ of length $k$ in the corresponding sequence of $G_0(f)$ and $G_1(f)$. 

Moreover, analogous statements hold if we replace top nodes with bottom nodes. Indeed, if there is 
a path from $a$ to $b$ of length $k$ in the infinite sequence of $G_0(f),G_1(f)$ then there is a path from $a_B$ to $b_B$ (in any of the two copies) of length $2k$ in the corresponding sequence of $g_0(f), g_1(f), g_2(f)$ where $g_2(f)$ appears at, and only at, the odd time steps. Similarly, a partial converse is that if there is a 
path from $a_B$, which is not one of the two $F$ nodes, to the $F$ node in its copy of length $2k$ in a sequence of $g_0(f), g_1(f), g_2(f)$ where $g_2(f)$ appears at, and only at, the odd time steps, then there is a path from $a$ to $F$ of length $k$ in the corresponding sequence of $G_0(f), G_1(f)$. 

Finally, we observe that in any sequence of $g_0(f), g_1(f), g_2(f)$, there is a path from any top node to the $F$ in its copy the first time $g_2(f)$ appears at an odd time step or one of $g_0(f), g_1(f)$ appears at an even time step. Similarly, there is a path from any bottom node to the $F$ in its copy the first
time $g_2(f)$ appears at an even time step or one of $g_0(f), g_1(f)$ appears at an odd time step. 

\smallskip

\noindent {\bf Proof that deciding whether the construction converges to consensus is NP-hard:} With these preliminary observations in place, we now proceed to the proof. Suppose $\{A_0(f), A_1(f)\}$ is not a consensus set. By Remark \ref{pathn1}, this means that 
there exists an infinite sequence of $G_0(f), G_1(f)$  such that there is no path from $S$ to $F$. Consider replacing every $G_0(f)$ by
$g_0(f)$, every $G_1(f)$ by $g_1(f)$ and inserting $g_2(f)$ into every even time slot. By the preceeding remarks,  we have that paths beginning at the 
two $S_T$ nodes never reach the $F$ nodes in their halves of the graph. Since the  only edges connecting the two halves of the 
graphs connect the $F$-nodes, we have that the set of reachable nodes from the two $S_T$ nodes never intersect.  By item $(5)$ of Lemma \ref{firstequivalence}, this implies $\{B_0(f), B_1(f), B_2(f)\}$ is not a consensus set. 

Conversely, suppose that $\{A_0(f), A_1(f)\}$ is a consensus set. We next show that for every pair of nodes $i,j$  and every sequence of $g_0(f), g_1(f), g_2(f)$ of length $4(n+1)+1$, there is some third node $k$ reachable by a path starting from both $i$ and $j$. This will prove that $\{B_0(f), B_1(f), B_2(f)\}$ is a consensus set, once again by item $(5)$ of Lemma \ref{firstequivalence}.

First observe that the conclusion true if one of $i,j$ is a top node and the other is a bottom node. Indeed, in this case the supports of the reachable
nodes from $i$ and $j$ intersect after three steps. 

Consequently, we may either assume both are top nodes or both are bottom nodes. Let us suppose both $i$ and $j$ are top nodes. Now consider any sequence of $g_0(f), g_1(f), g_2(f)$; it will either be true or false that, among the first $4(n+1)$ time steps,  $g_2(f)$ appears at, and only at, the even time steps.

Suppose it is true. Then there is a path from each $i,j$ to the $F$ node in its copy in $4(n+1)$ steps: 1) from any top node, at most $2(n+1)$ steps to reach either $F$ or the $S_T$ in its copy, and then $2(n+1)$ additional steps to go from $S_T$ to $F$ and, if necessary take self-loops at 
$F$ the remainder of the time. Note that here we used the assumption that $\{A_0(f), A_1(f)\}$ is not a consensus set, and consequently the $F$ node
must be reachable by a path of length $2(n+1)$ starting from $S_T$ regardless of the graph sequence. Consequently, the support of the reachable nodes from $i,j$ intersects after $4(n+1)+1$ steps. 

Suppose it is false. Then in the first $4(n+1)$ steps we have that either $g_2(f)$ appears at an odd time step or one of $g_0(f), g_1(f)$ appears at an even time
step. Since there is a path from every top node to $F$ node in its copy the first time either of these things happens, we obtain that the support of reachable nodes
from $i,j$ intersects in $4(n+1)+1$ steps in this case as well. 

The argument under the assumption both $i,j$ are both bottom nodes is similar, and this completes the proof. 


 {\vbox{\hrule height0.6pt\hbox{%
   \vrule height1.3ex width0.6pt\hskip0.8ex
   \vrule width0.6pt}\hrule height0.6pt
  }}

  \section{Conclusions\label{concl}} We have provided a combinatorial necessary and sufficient condition for consensus, the so-called ``avoiding sets condition." This theorem can be used to show that checking consensus of a finite set of matrices is NP-hard but decidable. This is the case even when there are just two matrices in the set. For undirected matrices, our NP-hardness result is for three or more matrices. 
  
Our work points to some intriguing open questions. We have left open the complexity of deciding whether a set of  two undirected matrices is a consensus set. Moreover, while the NP-hardness results conjecturally 
rule out a polynomial-time algorithm, it is not at present clear if it is 
possible to check whether a finite set of $n \times n$ matrices is a consensus set in time which is even singly exponential in $n$. Indeed, the naive algorithm based on the avoiding sets condition checks properties of 
all products of length $O(3^n)$, and the number of such products is doubly exponential in $n$. 

More broadly, even though stability determination for switched systems often turns out to be 
undecidable \cite{vb00}, it may turn out to be decidable (and even tractable) for many important 
subclasses of switched systems. It would be very interesting to relate the complexity of stability questions
such as the one we consider here to the structural properties of the systems in question.

\end{document}